\documentclass[12pt,a4paper,leqno]{amsart}
\usepackage[english]{babel}
\usepackage[T1]{fontenc}
\usepackage{amssymb}
\usepackage{amsfonts}
\usepackage{amsmath}
\usepackage{mathptm}

\swapnumbers
\title{A generalized mean value inequality for subharmonic functions and applications}
\author{Juhani Riihentaus}
\date{\noindent 21 February 2003}
\begin{document}
\pagestyle{myheadings}
\maketitle
\setcounter{page}{1}

\vspace*{0.1cm}

\noindent\begin{abstract}
\noindent{\footnotesize{
If $u\geq 0$ is subharmonic on a domain $\Omega$ in ${\mathbb R}^n$ and $p>0$,
then it is well-known
that there is a constant $C(n,p)\geq 1$ such that $u(x)^p\leq C(n,p)\, {\mathcal
MV}(u^p,B(x,r))$ for each ball $B(x,r)\subset \Omega$. We recently showed  that a similar result 
holds more generally for functions of the form $\psi \circ u$ where $\psi :{\mathbb R}_+\to
{\mathbb R}_+$ may be any surjective, concave function whose inverse $\psi ^{-1}$
satisfies the $\varDelta_2$-condition. Now we point out that this result can be extended slightly 
further. We also apply this extended result to the  weighted boundary behavior and nonintegrability
 questions 
of subharmonic and superharmonic functions.}}
\end{abstract}

\vspace{16pt} 

\footnotetext[1]{ 
(Secondary).}

\setcounter{equation}{0}

\section{Introduction}

\subsection*{1.1 Previous results}
If $u$ is a nonnegative and subharmonic function on $\Omega$ and $p>0$, then
there is a constant $C=C(n,p)\geq 1$ such that
\begin{equation}
u(x)^p\leq \frac{C}{m(B(x,r))}\int_{B(x,r)}u(y)^p\, dm(y)\end{equation}
for all $B(x,r)\subset \Omega$. Here $\Omega$ is a domain in ${\mathbb R}^n$, $n\geq 2$, $B(x,r)$ is the Euclidean
ball with center $x$ and radius $r$, and $m$ is the Lebesgue measure in ${\mathbb R}^n$.
See [FeSt72, Lemma~2, p. 172], [Ku74, Theorem~1, p. 529], [Ga81, Lemma~3.7, pp. 121-123],
[AhRu93, (1.5), p. 210]. These authors considered only the case when $u=\vert v\vert$ and 
$v$ is harmonic function.  However, the proofs in [FeSt72] and [Ga81] apply verbatim also in the general
 case of nonnegative subharmonic functions. This was pointed out in [Ri89, Lemma, p. 69],
[Su90, p. 271], [Su91, p. 113],
 [Ha92, Lemma~1, p. 113], [Pa94, p. 18] and [St98, Lemma 3, p. 305]. In [AhBr88, p. 132] it was pointed
out that a modification of the proof in [FeSt72] gives in fact a slightly more general result, see 
{\textbf{2.1}} below. A possibility for an essentially different proof was pointed out already in [To86, pp. 188-190]. Later other different proofs were given in [Pa94, p. 18, and Theorem 1, p. 19] (see also [Pa96, Theorem A, p. 15]), [Ri99, Lemma 2.1, p. 233] and [Ri01,Theorem, p. 188].
The results in [Pa94], [Ri99] and [Ri01] hold in fact 
for more general function classes than just for nonnegative subharmonic functions. See {\textbf{2.1}} and 
Theorem A  below.
 Compare also [DBTr84] and [Do88, p. 485]. 

The inequality $(1)$ has many applications. Among others, it has been
applied to
the (weighted) boundary behavior of nonnegative subharmonic functions
[To86, p. 191], [Ha92,
Theorems~1 and 2, pp. 117-118], [St98, Theorems 1, 2 and 3, pp. 301, 307],
[Ri99, Theorem, p. 233] and on the nonintegrability of
subharmonic and superharmonic functions
[Su90, Theorem
2, p.
271],
[Su91, Theorem,
 p.
113].

Because of the importance of the mean value inequality $(1)$, it is worthwhile to present a unified result 
which contains this mean value inequality and all its generalizations cited above. Below in  Theorem~2.5 we propose 
such a generalization. Instead of nonnegative
subharmonic functions we formulate our result slightly more generally for functions 
which we call \emph{quasi-nearly subharmonic functions} and which will be defined in {\textbf{2.1}}. 
We also give
two applications.

As the first application we improve in Theorem~3.4 below our recent result [Ri99, Theorem, p. 233] on the 
weighted boundary behavior of nonnegative
subharmonic functions. As the  second application we give in Corollary 4.5 
below  a supplement to Suzuki's results on the
nonintegrability of superharmonic and subharmonic functions [Su91, Theorem, p. 113]. Our result is  a limiting case 
to Suzuki's results.

\subsection*{1.2 Notation} Our notation is more or less standard, see
[Ri99]. However, for convenience of the reader we recall here the following.
We use the common convention $0\cdot \infty =0$. We write $\nu_n=m(B(0,1))$.
The
$d$-dimensional Hausdorff (outer) measure in ${\mathbb R}^n$ is denoted by $H^d$,
$0\leq d\leq n$. In the sequel $\Omega$ is always a domain in ${\mathbb R}^n$,
$\Omega \ne {\mathbb R}^n$, $n\geq 2$. The diameter of $\Omega$ is denoted by 
 diam\,$\Omega$.
The distance from $x\in \Omega$ to $\partial \Omega$, the boundary of $\Omega$,
is denoted by $\delta (x)$. ${\mathcal L}^1_{\mathrm{loc}}(\Omega )$ is the space of
locally integrable functions on $\Omega$. Our constants $C$ are always
positive, mostly $\geq 1$, and they may vary from line to line.

\section{Quasi-nearly subharmonic functions}
\subsection*{2.1 The definition}
  We call a
(Lebesgue) measurable function
$u:\Omega
\to
 [-\infty ,\infty ]$
 {\textit {quasi-nearly subharmonic}}, if $u\in {\mathcal L}^1_{{\textrm{loc}}}(\Omega
 )$ and if there is a constant $C_0=C_0(n,u,\Omega )\geq 1$ such that
\begin{equation}
u(x)\leq \frac{C_0}{r^n}\int_{B(x,r)}u(y)\, dm(y)\end{equation}
for any ball $B(x,r)\subset \Omega$. Compare [Ri99, p. 233]
and [Do57, p. 430]. Nonnegative quasi-nearly subharmonic
functions have previously been considered by [Pa94]
(he called them "functions satisfying the sh$_{\mathrm{K}}$-condition") and in [Ri99] 
(where they were called "pseudosubharmonic functions").
See [Do88, p. 485] for an even  more general function class
of (nonnegative) functions. As a matter of fact, also we will restrict ourselves
 to nonnegative functions.

Nearly
subharmonic functions, thus also quasisubharmonic and subharmonic functions, are examples
of quasi-nearly subharmonic functions.
Recall that a
function $u\in {\mathcal L}^1_{\textrm{loc}}(\Omega )$ is nearly subharmonic, if $u$
satisfies $(2)$ with $C_0=\frac{1}{\nu_n}$.
See [Her71, pp. 14, 26]. Furthermore, if $u\geq 0$ is
subharmonic and $p>0$, then by $(1)$ above, $u^p$ is quasi-nearly subharmonic.
By [Pa94, Theorem~1, p. 19] or [Ri99, Lemma 2.1, p. 233] this holds even if $u\geq 0$
is quasi-nearly subharmonic. See also [AhBr88, p. 132].

\subsection*{2.2 Permissible functions}
In
[Ri01, Theorem, p. 188]
we proved the following result which contains essentially the cited results in [Pa94] and [Ri99] 
 as special cases.

\vspace{4pt}

\noindent\textbf{Theorem A.} ([Ri01, Theorem, p. 188]){\emph{ Let
$u$ be a  nonnegative subharmonic function on $\Omega$. Let  $\psi :{\mathbb R}_+\to 
{\mathbb R}_+$ be a concave surjection whose inverse $\psi^{-1}$  satisfies the
$\varDelta_2$-condition. Then there exists a constant $C=C(n,\psi ,u)\geq 1$
such that
\[ \psi (u(x_0))\leq \frac{C}{{\varrho}^n}\int_{B(x_0,\varrho )}\psi (u(y))\,
dm(y)\] for any ball $B(x_0,\varrho )\subset \Omega$.}}

\vspace{4pt}

Recall that a function $\psi :{\mathbb R}_+\to {\mathbb R}_+$ satisfies the
$\varDelta_2$-{\textit {condition}}, if there is a constant $C=C(\psi )\geq 1$
such that
\begin{displaymath}\psi (2t)\leq C\, \psi (t)\end{displaymath} for all $t\in {\mathbb R}_+$.

In order to improve our result, Theorem A above, still further,
we give the following definition.
A function $\psi :{\mathbb {R}}_+\to {\mathbb {R}}_+$ is {\textit {permissible}},
if
there  is a nondecreasing, convex function
$\psi_1
:{\mathbb
{R}}_+\to
{\mathbb
{R}}_+$
and an increasing surjection
$\psi_2
:{\mathbb
{R}}_+\to
{\mathbb
{R}}_+$ such that $\psi =\psi_2\circ \psi_1$ and such that the following
conditions are satisfied:
\begin{itemize}
\item[(a)] $\psi_1$ satisfies the
$\varDelta_2$-condition.
\item[(b)] $\psi_2^{-1}$ satisfies the
$\varDelta_2$-condition.
\item[(c)] The function $t\mapsto \frac{t}{\psi_2(t)}$ is
{\textit {quasi-increasing}}, i.e. there is a constant $C=C(\psi_2)\geq 1$ such
that
\begin{displaymath}\frac{s}{\psi_2(s)}\leq C\,    \frac{t}{\psi_2(t)}\end{displaymath}
for all $s, t\in {\mathbb R}_+$, $0\leq s\leq t$.
\end{itemize}
Observe that the condition (b) is equivalent with the following condition.
\begin{itemize}
\item[(b')] For some constant $C=C(\psi_2)\geq 1$,
\begin{displaymath}\psi_2(Ct)\geq 2\, \psi_2(t)\end{displaymath}
for all $t\in {\mathbb R}_+$.
\end{itemize}
If $\psi$ is a permissible function,  we will in the sequel use always
one, fixed constant
 $C_1=C_1(\psi )$, in the case of  all the properties (a),
(b),
(c)
(and
(b')).
If $\psi : {\mathbb
R}_+\to {\mathbb R}_+$ is an increasing surjection satisfying the conditions (b)
and (c), we say that it is {\textit{strictly permissible}}. Permissible functions are necessarily continuous.

\vspace{4pt}

Let it be noted that the condition (c) above is indeed 
natural. For just one counterpart
to it,  see e.g.  [HiPh57, Theorem 7.2.4, p. 239].

Observe that our previous definition for permissible functions in [Ri99, 1.3, p. 232] was {\emph{much more restrictive}}:
 A function $\psi :{\mathbb{R}}_+\rightarrow {\mathbb{R}}_+$ was there 
 defined to be {\emph{permissible}} if it is of the form $\psi (t)=\vartheta (t)^p$,
 $p>0$, where $\vartheta :{\mathbb{R}}_+\rightarrow {\mathbb{R}}_+$ is a nondecreasing, convex
function satisfying the $\varDelta_2$-condition.

\vspace{4pt}

{\textbf{2.3~Remarks}} The following list gives  examples of permissible functions (we leave the slightly tedious 
verifications to the reader). Our list, especially (ii),
 shows that the considered class of permissible functions is wide and
natural. On the other hand, in view of the simple example in (vi), one sees that functions of type (ii)
 are by no means the only permissible functions: There exists a huge amount of
permissible functions of other types. 
   
\begin{itemize}
\item[(i)]  The functions $\psi_1(t)=\vartheta (t)^p$, $p>0$, already  considered
in [Ri99, (1.3), p. 232, and Lemma 2.1, p. 233].
\item[(ii)]  Functions of the form $\psi_2=\phi_2 \circ \varphi_2$,
where
$\phi_2:{\mathbb R}_+\to {\mathbb R}_+$ is a concave surjective function whose
inverse $\phi_2^{-1}$ satisfies the $\varDelta_2$-condition, and
$\varphi_2:{\mathbb R}_+\to {\mathbb R}_+$ is a nondecreasing convex
function
 satisfying the
$\varDelta_2$-condition. (Observe here that any concave function
$\phi_2:{\mathbb R}_+\to {\mathbb R}_+$ is necessarily
nondecreasing.) These (or, to be more exact, the functions $\phi_2$ defined above) were considered
in [Ri01, Theorem, p. 188].
\item[(iii)]
 $\psi_3(t)=c\, t^{p\,\alpha}[\log (\delta +t^{p\,\gamma})]^{\beta}$, where
$c>0$, $0<\alpha <1$, $\delta \geq 1$, and $\beta ,\gamma \in {\mathbb {R}}$
are such that $0<\alpha +\beta \gamma <1$, and $p\geq 1$.
\item[(iv)]  For $0<\alpha <1$, $\beta \geq 0$ and $p\geq 1$,
\begin{displaymath}\psi_4(t)=\left\{ \begin{array}{ll}p^{\beta}t^{p\, \alpha},&{\textrm{ for }}
 0\leq t\leq e,\\
t^{p\, \alpha}(\log t^p)^{\beta},&{\textrm{ for }} t>e.\end{array}\right. \end{displaymath}
\item[(v)] For $0<\alpha <1$, $\beta < 0$ and $p\geq 1$,
\begin{displaymath}\psi_5(t)= \left \{ \begin{array}{ll}
 (\frac{-\beta \,p}{\alpha})^{\beta} t^{p\,\alpha},
&{\textrm{ for }} 0\leq t\leq e^{-{\beta} /{\alpha}},\\
 t^{p\,\alpha}(\log t^p)^{\beta},
&{\text{ for }}
 t>e^{-\beta /{\alpha}}.\end{array} \right. \end{displaymath}
\item[(vi)] For  $p\geq 1$,
\begin{displaymath}\psi_6(t)= \left \{ \begin{array}{lll}
 2n+\sqrt{t^p-2n},&{\textrm{ for }} t^p\in [2n,2n+1), &n=0,1,2,\ldots,\\
 2n+1+[t^p-(2n+1)]^2,
&{\text{ for }}
 t^p\in [2n+1,2n+2), &n=0,1,2,\ldots .\end{array} \right. \end{displaymath}
\end{itemize}

For $p=1$ the functions in (i), (iii), (iv), (v), (vi), and also in (ii) provided $\varphi_2(t)=t$,
 are 
strictly permissible.

\vspace{4pt}

\noindent{\textbf{2.4 Remark.}} Our previous results were restricted  to the cases 
where $\psi$ was either of \mbox{type (i)} ([Ri99, (1.3), p. 232, and Lemma 2.1, p. 233])  or of
 type (ii) ([Ri01, Theorem, p. 188] (or Theorem A above)). Now we give a refinement to our results
in a unified form. The proof is a
modification of Pavlovi\'c's argument [Pa94, proof of
Theorem~1, p.
20].

\vspace{4pt}

\noindent\textbf{2.5 Theorem.} {\emph{Let $u$ be a  nonnegative quasi-nearly subharmonic function on $\Omega$. 
If $\psi :{\mathbb R}_+\to {\mathbb R}_+$ is a permissible function, then $\psi \circ u$ is
quasi-nearly subharmonic on $\Omega$.}}

\vspace{4pt}

{\textit {Proof}}. In view of [Ri99, Lemma 2.1, p. 233] we may restrict us to the case where
$\psi =\psi_2:{\mathbb R}_+\to {\mathbb R}_+$ is strictly permissible.

Since $\psi$ is  continuous, $\psi \circ u$ is measurable and 
$\psi \circ u\in  {\mathcal{L}}^1_{\textrm{loc}}(\Omega)$. It remains to show that $\psi \circ u$ satisfies 
the generalized mean value inequality (2).
But this can be seen exactly as in [Ri01, proof of Theorem, pp. 188-189], the only 
difference being that  instead of the property 2.4 in 
[Ri01, p. 188] of concave functions, one now  uses  the above property (c) in {\textbf{2.2}} 
of permissible functions.
\null{} \quad \hfill $\square$

\vspace{4pt}

\section{Weighted boundary behavior}
\subsection*{3.1 Stoll's result.} Improving previous results of Gehring [Ge57, 
Theorem~1, p. 77] and
Hallenbeck [Ha92, Theorems~1 and 2, pp. 117-118], Stoll [St98,
Theorem~2, p. 307] gave the following result.

\vspace{4pt}

\noindent\textbf{Theorem B.}
{\emph{ Let $f$ be a nonnegative subharmonic function on a domain  $G$ in ${\mathbb R}^n$,
$G\ne {\mathbb R}^n$, $n\geq 2$, with ${\mathcal C}^1$ boundary. Let
\begin{equation} \int_{G}f(x)^p\delta (x)^{\gamma}\, dm(x)<\infty \end{equation}
for some $p>0$ and $\gamma >-1-\beta (p)$. Let $0<d\leq n-1$. Then for each
$\tau
\geq
1$ and $\alpha >0$ ($\alpha >1$ when $\tau =1$), there exists a subset
$E_{\tau}$ of
$\partial
G$ with
$H^d(E_{\tau})=0$ such that
\begin{equation*} \lim_{\rho\to 0}\{ \sup_{x\in \Gamma_{\tau ,\alpha ,\rho}(\zeta )}[\delta
(x)^{n+\gamma -\frac{d}{\tau}}f(x)^p]\} =0\end{equation*}
for all $\zeta \in \partial G \setminus E_{\tau}$.}}

\vspace{4pt}

Above, for $\zeta \in \partial G$ and
$\varrho >0$,
\[ \Gamma_{\tau ,\alpha ,\varrho}(\zeta )=
     \Gamma_{\tau ,\alpha }(\zeta )\cap G_{\varrho}, \]
where
\[ \Gamma_{\tau ,\alpha}(\zeta )=\{\, x\in G:\, \vert x- \zeta \vert ^{\tau}
<\alpha \, \delta
(x)
\,\},\quad
G_{\varrho}= \{\, x\in G:\, \delta (x)<\varrho \,\}. \]
Moreover,
$\beta (p)=\max \{\, (n-1)(1-p),0\, \}$. Stoll makes the assumption $\gamma
>-1-\beta (p)$ in order to exclude the trivial case $f\equiv 0$. As a matter of
fact, it follows from a result of Suzuki [Su90, Theorem~2, p.
271] that $(3)$ together with the condition $\gamma \leq -1-\beta (p)$ implies
indeed that $f\equiv 0$, provided $G$ is a bounded domain with $\mathcal{C}^2$ 
boundary.

Our previous improvement 
[Ri99, Theorem, p. 233] to Stoll's result, Theorem B above, can now be refined slightly further, just using our 
above result Theorem 2.5. This refinement will be given in Theorem 3.4 below.
For this purpose we first recall, and just for the
convenience of the reader, some terminology from [Ri99, pp.
231--232].
\subsection*{3.2 Admissible functions} A function
$\varphi :{\mathbb {R}}_+\to {\mathbb {R}}_+$ is {\it {admissible}}, if
it is increasing (strictly), surjective, and  there are constants
$C_2>1$ and
$r_2>0$ such that
\begin{equation} \varphi
(2t)\leq C_2\, \varphi
(t)
\, \, \, \, 
{\textrm{ {and}}}\, \, \, \, \, \,  \varphi ^{-1} (2s)\leq C_2\, \varphi ^{-1}(s)
\, \, \, \, {\textrm  {for all}}
\, \, \, \, \, \, s, \, t\in {\mathbb {R}}_+,\, \, 0\leq s,\, t\leq
r_2. \end{equation}

Nonnegative,  nondecreasing functions $\varphi_1(t)$ which satisfy
 the $\varDelta_2$-condition and for which the functions
$t\mapsto \frac{\varphi_1(t)}{t}$ are nondecreasing,
are examples of admissible functions.
Further examples are
$\varphi_2(t)=c\, t^{\alpha}[\log (\delta +t^{\gamma})]^{\beta},$
where
$c>0$,
$\alpha
>0$,
$\delta
\geq
1$,  and
$\beta
,\gamma
\in
{\mathbb
{R}}$  are such that
$\alpha +\beta \gamma >0$.

\subsection*{3.3 Accessible boundary points  and approach
regions} Let $\varphi :{\mathbb {R}}_+\to {\mathbb {R}}_+$ be an
admissible function and let $\alpha >0$. We say that $\zeta \in \partial
\Omega$ is
$(\varphi ,\alpha )$-{\it {accessible}}, if
\[ \Gamma_{\varphi}(\zeta ,\alpha )\cap B(\zeta ,\rho )\ne
\emptyset \] for all
$\rho
 >0$. Here
\[ \Gamma_{\varphi} (\zeta ,\alpha )=\{\, x\in \Omega \, :\, \varphi (\vert
x-\zeta \vert )<\alpha \, \delta (x)\, \},\]
and we call it a  $(\varphi ,\alpha )$-{\textit {approach region
in}}
$\Omega$ {\textit {at}} $\zeta$.

Mizuta [Mi91] has considered boundary limits of harmonic functions in
Sobolev-Orlicz classes on bounded Lipschitz domains $U$ of ${\mathbb R}^n$, $n\geq
2$. His approach regions are of the form
\[ \Gamma_{\phi} (\zeta ,\alpha )=\{\, x\in U \, :\, \phi (\vert
x-\zeta \vert )<\alpha \, \delta (x)\, \},\]
where now $\phi :{\mathbb R}_+\to {\mathbb R}_+$ is a nondecreasing function which
satisfies the $\varDelta_2$-condition and is such that $t\mapsto \frac{\phi
(t)}{t}$ is nondecreasing. As pointed out above, such functions are admissible
in our sense. In fact, they form a proper subclass of our admissible functions.

\vspace{4pt}

\noindent\textbf{3.4 Theorem.}
{\emph{Let $H^d(\partial
\Omega )<\infty$ where $0\leq d\leq n$. Suppose that
$u$ is a nonnegative quasi-nearly subharmonic function in $\Omega$.  Let
$\varphi
:{\mathbb {R}}_+\to {\mathbb {R}}_+$ be an admissible function and $\alpha >0$. Let
$\psi
:{\mathbb
{R}}_+\to {\mathbb {R}}_+$ be a permissible function. Suppose that
\begin{equation} \int_{\Omega}\psi (u(x))\delta (x)^{\gamma}\, dm(x)<\infty \end{equation}
for some $\gamma \in {\mathbb {R}}$.
Then
\[ \lim_{\rho\to 0}( \sup_{x\in \Gamma_{\varphi ,\rho}(\zeta , \alpha )}\{
\delta
(x)^{n+\gamma}\,
[\varphi
^{-1}(\delta
(x))]^{-d}\,
\psi
(u(x))\}
)
=0\] for $H^d$-almost every $(\varphi ,\alpha )$-accessible point  $\zeta \in
\partial
\Omega$. Here
\[ \Gamma_{\varphi ,\rho}(\zeta ,\alpha )=\{\, x\in \Gamma_{\varphi}(\zeta
,\alpha )\, :\, \delta (x) <\rho \, \}.\]}}

The proof is verbatim the same as [Ri99, proof of
Theorem, pp.
235--238], except that now we just replace [Ri99, Lemma~2.1,
p. 233] 
by the more general Theorem~2.5 above.
\null{} \quad \hfill \qed

\vspace{4pt}

\noindent {\textbf{3.5 Remark.}} Unlike Stoll, we have imposed no restrictions on the
exponent $\gamma$ in order to exclude the trivial case $u\equiv 0$. We
refer, however, to such possibilities in Remark 4.7 below, after having given in Corollary 4.5 a limiting case result 
 for a result of Suzuki.

\section{A Limiting case result  to nonintegrability results of Suzuki}

\subsection*{4.1 Suzuki's result} Suzuki  [Su91, Theorem and 
its proof, pp. 113--115] gave the following result.

\vspace{4pt}

\noindent\textbf{Theorem~C.}
{\emph{Let $0<p\leq 1$. If a superharmonic (resp. nonnegative subharmonic) function
$v$ on  $\Omega$ satisfies
\begin{equation} \int_{\Omega}\vert v(x)\vert ^p\, \delta (x)^{np-n-2p}\,
 dm(x)<\infty ,\end{equation}
then $v$ vanishes identically.
}}             

\vspace{4pt}

Suzuki pointed also  out  that his result is sharp in the following sense: If $p$, 
$0<p\leq 1$,
is fixed, then the exponent $\gamma =np-n-2p$ cannot be increased. On the other
hand, clearly $-n<\gamma \leq -2$, when $0<p\leq 1$. Since the
class of permissible
functions include, in addition  the functions $t^p$, $0<p\leq 1$, also a
large
amount of essentially different functions (see {\textbf{2.3}} above),
one is tempted to ask whether there exists any limiting case result for Suzuki's results, corresponding to 
the case $p=0$. 
To be more precise, we pose the following  question: 

{\emph{Let $\Omega$ and $v$ be as above. Let $\gamma \leq -n$
and let $\psi :{\mathbb{R}}_+\rightarrow {\mathbb{R}}_+$ be permissible. Does the condition
\[ \int_{\Omega}\psi (\vert v(x)\vert )\delta (x)^{\gamma}\, dm(x)<\infty ,\]
imply $v\equiv 0$?}}

Observe that the least severe form of  above integrability condition occurs when  $\gamma =-n$. 

Below in 
Corollary 4.5 we answer the question in the affirmative, in the case of any strictly permissible function $\psi$.
In order to achieve this, we first formulate below in Theorem 4.3 a general
result for arbitrary $\gamma \leq -2$ which is, for $-n< \gamma \leq -2$, however, essentially more or less just
 Suzuki's above result (see Remarks 4.4 (b) below).
Our formulation has the advantage that, unlike Suzuki's result, it contains  a certain limiting case, Corollary 4.5, too.

The
proof which we below write down (and in quite detail, just for the convenience
of the
reader) is merely a slight modification of Suzuki's argument, combined
with our version for the generalized mean value inequality (Theorem~2.5 above),
and also some additional estimates.

\vspace{4pt}

\noindent\textbf{4.2~Lemma.}
{\emph{Let $u$ be  a nonnegative subharmonic function on $\Omega$.
Suppose $\psi :{\mathbb R}_+\to {\mathbb R}_+$ is a permissible function such that
\[ \int_{\Omega}\psi (u(x))\delta (x)^{\gamma}\, dm(x)<\infty \]
for some $\gamma \in {\mathbb R}$. Then $\psi (u(x))=o(\delta (x)^{-n-\gamma})$ as
$\delta (x)\to 0$.}}

\vspace{4pt}

{\textit{Proof}}. By Theorem~2.5, $\psi \circ u$ is quasi-nearly subharmonic on
$\Omega$. Write for $x\in \Omega$, $B=B(x,\delta (x))$ and
$B_0=B(x,\frac{\delta (x)}{2})$. Since
\begin{equation} \frac{1}{2}\delta (x)<\delta (y)<\frac{3}{2}\delta (x)\end{equation}
for all $y\in B_0$, we get
\[ \begin{array}{rl} \psi (u(x))&\leq 2^{n+\vert \gamma \vert}\,C_0\, \delta
(x)^{-n-\gamma}\int_{B_0}\psi (u(y))\delta (y)^{\gamma}\, dm(y)\\
&\leq C\, \delta
(x)^{-n-\gamma}\int_{\Omega_{\delta (x)}}\psi (u(y))\delta (y)^{\gamma}\,
dm(y),\end{array} \]
where $C=C(\gamma ,n,\psi ,u)>0$ and \,$\Omega_{\delta (x)}=\{y\in \Omega \,:\, \delta (y)<\delta (x)\,\}$. The claim follows.
\null{
}\quad \hfill \qed

\vspace{4pt}

\noindent\textbf{4.3~Theorem.}
{\emph{Let $\Omega$ be bounded. Let $v$ be
a superharmonic (resp. nonnegative subharmonic) function  on $\Omega$. Let
$\psi :{\mathbb R}_+\to {\mathbb R}_+$ be a strictly permissible function. Suppose
\begin{equation} \int_{\Omega}\psi (\vert v(x)\vert )\delta (x)^{\gamma}\, dm(x)< \infty ,\end{equation}
 where $\gamma\leq -2$ is such that there is a constant
$C=C(\gamma ,n,\psi ,\Omega )>0$
 such that}}
\begin{equation} s^{n+\gamma}\leq \psi (C\, s^{n-2}) \,{\textit{ for all }}\,
s>\frac{1}{{\textrm{diam}}\,\Omega }.\end{equation} {\emph{Then $v$ vanishes identically.}}

\vspace{4pt}

\noindent{\textbf{4.4~Remarks.}} Next we consider the assumptions in Theorem 4.3. 
\begin{itemize}
\item[(a)]  Our assumption $\gamma \leq -2$ is unnecessary, and it could be dropped:
 If
$\gamma \in {\mathbb R}$, then  it follows easily  from
$(9)$ and from the property
(c) in {\textbf{2.2}}
of  strictly permissible functions that indeed $\gamma \leq -2$.
\item[(b)] Suppose that $-n<\gamma \leq -2$. If, instead of $(9)$, one supposes that
\begin{equation*} s^{n+\gamma}\leq \psi (C\, s^{n-2}) {\textit{ for all }}
s>0,\end{equation*}
then  clearly
\[ \psi(\vert v(x)\vert )\geq C^{-\frac{n+\gamma}{n-2}}\vert v(x)\vert ^{\frac{n+\gamma}{n-2}}\]
for all $x\in \Omega$. Thus $(8)$ implies that 
\[ \int_{\Omega}\vert v(x)\vert ^{\frac{n+\gamma}{n-2}}\delta (x)^{\gamma}\, dm(x)<\infty ,\]
and hence  $v \equiv 0$ by Suzuki's result, Theorem C above. Recall that here $0<p=\frac{n+\gamma}{n-2}\leq 1$ and $\gamma =np-n-2p$.
Thus Theorem 4.3, but now the assumption $(9)$ replaced with the aforesaid assumption, is just a restatement of 
Suzuki's result for bounded domains.
\item[(c)] If $\gamma \leq -n$, then the condition $(9)$ clearly holds, since $\psi$ is
strictly permissible. This case gives indeed  the already referred limiting  case for Suzuki's result. See Corollary 4.5 below.
\end{itemize}

\vspace{4pt}

{\textit{Proof of Theorem~4.3}}.   We write the proof down only for the case $n\geq
3$. Write $v^+=\max \{v,0\}$ and $s=v^-=-\min \{v,0\}$. Then
$\vert v\vert
=v^++v^-$ and
$s\geq 0$ is subharmonic. (Resp. if $v$ is nonnegative and subharmonic, let
$s=v$.) Proceeding as Suzuki, but using also some additional estimates, we
will show that
$s\equiv
0$.

By $(8)$,
\[ \int_{\Omega}\psi (s(x))\delta (x)^{\gamma}\,dm(x)<\infty ,\]
thus
\[ v_{\psi ,\gamma }(x)=\int_{\Omega}G_{\Omega}(x,y)\psi (s(y))\delta
 (y)^{\gamma}\, dm(y)< \infty \]
is a potential by [Hel69, Theorem~6.3, p. 99]. Here $G_{\Omega }$ is the
Green function of $\Omega$.

By Lemma~4.2, $\psi (s(y))\leq C\, \delta (y)^{-n-\gamma}$. Thus also

\begin{equation} s(y)\leq \psi ^{-1}(C\, \delta (y)^{-n-\gamma})\leq C\, \psi^{-1}(\delta
(y)^{-n-\gamma})\end{equation}
for all $y\in \Omega$, where $C=C(\gamma ,n,\psi ,s,\Omega )\geq 1$. Let $x\in
\Omega$ be fixed for a while. Let $B=B(x,\delta (x))$ and
$B_0=B(x,\frac{\delta
(x)}{2})$. Using $(10)$ (and $(7)$) one gets
\[ s(y)\leq C\, \psi ^{-1}(\delta (y)^{-n-\gamma})\leq C\,
\psi^{-1}(2^{n+\vert \gamma \vert}\delta
(x)^{-n-\gamma})\leq C\, \psi^{-1}(\delta (x)^{-n-\gamma})\]
for all $y\in B_0$. Therefore
\[ \frac{\psi (s(y))}{s(y)}\geq C_1\,\frac{\psi (C\,\psi^{-1}(\delta
 (x)^{-n-\gamma}))}{C\,\psi^{-1}(\delta (x)^{-n-\gamma})}\geq C\, \frac{\delta
 (x)^{-n-\gamma}}{\psi^{-1}(\delta (x)^{-n-\gamma})}\]
for all $y\in B_0$. With the aid of this and of  a standard estimate for
the Green function $G_B(x,\cdot )$ in $B_0$ and of $(10)$, one gets
\[ \begin{array}{rl} v_{\psi ,\gamma }(x)&=\int_{\Omega}G_{\Omega}(x,y)\psi (s(y))\delta
 (y)^{\gamma}\,dm(y)
\geq \int_{B_0}G_{B}(x,y)\psi (s(y))\delta
 (y)^{\gamma}\,dm(y)\\
&\geq \int_{B_0}G_{B}(x,y)s(y)\,\frac{\psi (s(y))}{s(y)}\delta
(y)^{\gamma}\,dm(y)\\
&\geq C\, \frac{\delta
(x)^{-n-\gamma}}{\psi^{-1}(\delta
(x)^{-n-\gamma})}\,\int_{B_0}\vert x-y\vert ^{2-n}
\delta
 (y)^{\gamma}s(y)\,dm(y)\\
&\geq C\, \frac{1}
{\delta (x)^{n-2}\psi^{-1}(\delta
 (x)^{-n-\gamma})}\, s(x).\end{array}\]
By $(9)$ we see that there is a constant $C_3\geq 1$ such that
\[ \delta(y)^{n-2}\,\psi^{-1}(\delta (y)^{-n-\gamma})\leq C_3\]
for all $y\in \Omega$. Combining this with the above estimate for $v_{\psi
,\gamma}$, one gets
\[ v_{\psi ,\gamma}(x)\geq C\, s(x),\]
where $C=C(\gamma ,n,\psi ,s,\Omega )>0$. Remembering that $x\in \Omega$ was
arbitrary, that $v_{\psi ,\gamma}$ is a potential and $s$ subharmonic, it
follows
from [Hel69, Corollary~6.19, p. 117] that $s\equiv 0$. Thus $v=v^+\geq 0$.
It remains to show that $v\equiv 0$.

As above,
\[ \begin{array}{rl}\int_{\Omega}G_{\Omega}(x,y)\, {\delta (y)}^{-2}\, dm(y)&\geq
\int_{B_0}G_B(x,y)\, {\delta (y)}^{-2}\, dm(y)\\
&\geq C\,\int_{B_0}{\vert x-y\vert }^{2-n}\, {\delta (y)}^{-2}\,dm(y)\\
&\geq C \,{\delta (x)}^{2-n}\, {\delta (x)}^{-2}\,{\nu}_n \,
(\frac{\delta (x)}{2})^n\, =\, C,
\end{array}\]
where $C=C(n)>0$. Thus by [Hel69, Lemma 6.1, p. 98, and Corollary 6.19, p. 117],
\begin{equation} \int_{\Omega}G_{\Omega}(x,y)\,\delta (x)^{-2}\, dm(y)=\infty \end{equation}
for all $x\in \Omega$. Consider next an arbitrary potential $w$ on $\Omega$,
\begin{equation} w(x)=\int_{\Omega}G_{\Omega}(x,y)\,d\lambda (y)\end{equation}
where $\lambda \ne 0$ is a measure on $\Omega$. From $(11)$ it follows that
\[ \begin{array}{rl} \int_{\Omega}w(x)\, \delta (x)^{-2}\, dm(x)&=
\int_{\Omega}[\int_{\Omega}G_{\Omega}(x,y)\,d\lambda (y)] \,\delta
(x)^{-2}\,dm(x)\\
&=\int_{\Omega}[\int_{\Omega}G_{\Omega}(x,y)\,\delta (x)^{-2}dm(x)]\,d\lambda
(y)=\infty.\end{array}\]
Using this and the facts that $\Omega$ is bounded and $w$, as a superharmonic
function, is locally integrable, one sees that
\begin{equation} \int_{\Omega_1}w(x)\, \delta (x)^{-2}\, dm(x)=\infty \end{equation}
where $\Omega_1 =\{\, x\in \Omega :\, \delta (x)<1\,\}.$

Suppose in particular that $w$ in $(12)$ is the potential of  the superharmonic
function $v_M=\inf \{ v,M\}$, where $M>0$. Then $v_M\geq w$, and one has by
$(13)$ and by the fact that $\gamma \leq -2$, 
\[ \begin{array}{rl} \infty >\int_{\Omega}\psi (v(x))\,\delta (x)^{\gamma}\, dm(x)&
\geq \int_{\Omega}\psi (v_M(x))\,\delta (x)^{\gamma}\,dm(x)\\
&\geq \int_{\Omega}v_M(x)\frac{\psi (v_M(x))}{v_M(x)}\,\delta
(x)^{\gamma}\,dm(x)\\
&\geq C\, \frac{\psi (M)}{M}\int_{\Omega}v_ M(x)\,\delta
(x)^{\gamma}\,dm(x)\\
&\geq C\, \frac{\psi (M)}{M}\int_{\Omega_1}w(x)\,\delta
(x)^{-2}\,dm(x)=\infty ,\end{array}\]
a contradiction unless $w\equiv 0$. Since $w\equiv 0$, the nonnegative superharmonic functions $v_M$, $M>0$, are 
in fact harmonic, e.g. by the Riesz Decomposition Theorem [Hel69, Theorem 6.18, p. 116]. Since this is impossible, one has 
$v\equiv 0$,
concluding the proof. 
\null{
}\quad \hfill \qed

\vspace{4pt}

\noindent\textbf{4.5~Corollary.}
{\emph{Let $\Omega$ be bounded. Let $v$ be
a superharmonic (resp. nonnegative subharmonic) function  on $\Omega$. Let
$\psi :{\mathbb R}_+\to {\mathbb R}_+$ be  any strictly permissible function and
let $\gamma \leq -n$. If
\begin{equation*} \int_{\Omega}\psi (\vert v(x)\vert )\delta (x)^{\gamma}\, dm(x)<\infty
,\end{equation*} then $v$ vanishes identically.}}

\vspace{4pt}

For the proof observe that the condition $(9)$ is indeed satisfied for $\gamma \leq
-n$, since $\Omega$ is bounded and $\psi$ is increasing.
\null{
}\quad \hfill \qed

\vspace{4pt}

\noindent {\textbf{4.6~Remark.}} The result of Theorem~4.3 does not, of course,  hold any
more,
if one replaces strictly permissible functions by permissible functions. For a
counterexample, set, say,  $v(x)=\vert x\vert ^{2-n}$, $\psi (t)=t^p$,
 where $\frac{n-1}{n-2}<p<\frac{n}{n-2}$, $\gamma =np-n-2p$ or just $\gamma
>1$. Then clearly
\[ \int_Bv(x)^p\,\delta (x)^{\gamma}\, dm(x)<\infty \]
but $v\not\equiv 0$.

\vspace{4pt}

\noindent {\textbf{4.7~Remark.}} Provided $\Omega$ is bounded and $\psi$ is strictly
permissible, one can, with the aid of Theorem 4.3 and Corollary 4.5,  exclude some trivial cases $u\equiv 0$ from the result of
Theorem~3.4  by imposing certain restrictions on the exponent $\gamma$.
 We point out only two cases:
\begin{itemize}
\item[(i)] By Corollary~4.5, $\gamma >-n$, regardless of $\psi$.
\nopagebreak[4]
\item[(ii)]  By Suzuki's result, Theorem~C above, $\gamma >np-n-2p$, in
the case when $\psi (t)=t^p$, $0<p\leq 1$.
\end{itemize}
\vskip1cm
\bibliographystyle{alpha}
\begin{center}\textsc{References}\end{center}

\vspace{8pt}

\footnotesize{
\begin{enumerate}

\item[{[AhBr88]}] P. Ahern and J. Bruna, \emph{Maximal and area integral characterizations of Hardy--Sobolev
spaces in the unit ball of} ${\mathbb{C}}^n$, Revista Matemática Iberoamericana \textbf{4}(1988), 123--153.
\item[{[AhRu93]}] P. Ahern and W. Rudin, \emph{Zero sets of functions in harmonic 
Hardy spaces}, Math. Scand.  \textbf{73}(1993), 209--214.
\item[{[DBTr84]}] E. Di Benedetto and N.S. Trudinger, \emph{Harnack inequalities for quasi-minima of variational 
integrals}, Ann. Inst. H. Poincaré, Analyse Nonlineaire  \textbf{1}(1984), 295--308.
\item[{[Do57]}] Y. Domar,\emph{ On the existence of a largest subharmonic minorant of a given
function}, Arkiv f\"or Mat. \textbf{3}(1957), 429--440.
\item[{[Do88]}] Y. Domar, \emph{Uniform boundedness in families related to subharmonic
 functions},  J. London Math. Soc. (2) \textbf{38}(1988), 485--491.
\item[{[FeSt72]}] C. Fefferman and E.M. Stein, \emph{H$^p$ spaces of several variables},
Acta Math.  \textbf{l29}(1972), 137--193.
\item[{[Ga81]}] J.B. Garnett, \emph{Bounded analytic functions}, Academic Press, 
 New York, 1981.
\item[{[Ge57]}]  F.W. Gehring, \emph{On the radial order of subharmonic functions}, 
J.\ Math. Soc. Japan \textbf{9}(1957), 77--79.
\item[{[Ha92]}] D.J. Hallenbeck, \emph{Radial growth of subharmonic functions}, 
Pitman Research Notes 262, 1992, pp. 113--121.
\item[{[Hel69]}] L.L. Helms, \emph{Introduction to potential theory}, Wiley-Interscience,
New York, 1969.
\item[{[Her71]}] M. Hervé, \emph{Analytic and Plurisubharmonic Functions in Finite and Infinite
Dimensional Spaces},  Lecture Notes in Mathematics 198,  Springer-Verlag, Berlin, 1971.
\item[{[HiPh57]}] E.~Hille and R.S.~Phillips, \emph{Functional Analysis and Semigroups}, 
American Mathematical Society, Colloquium publications XXXI, Providence, R.I., 1957.
\item[{[Ku74]}] \"U. Kuran, \emph{Subharmonic behavior of 
$\vert h\vert ^p$ ($p>0$, $h$ harmonic)},  J. London Math. Soc. (2)  \textbf{8}(1974),
529--538.
\item[{[Mi91]}] Y. Mizuta, \emph{Boundary limits of harmonic functions in Sobolev-Orlicz
classes},  Potential Theory (ed. M. Kishi), Walter de Gruyter\&Co, 
 Berlin $\cdot$ New York, 1991, pp. 235--249.
\item[{[Pa94]}] M. Pavlovi\'c, \emph{On subharmonic behavior and oscillation of functions
 on balls in} ${\mathbb R}^n$, 
 Publ. de l'Institut Mathém., Nouv. sér.  \textbf{55(69)}(1994), 18--22.
\item[{[Pa96]}] M. Pavlovi\'c, \emph{Subharmonic behavior of smooth functions},
 Math. Vestnik  \textbf{48}(1996), 15--21.
\item[{[Ri89]}] J. Riihentaus, \emph{On a theorem of Avanissian--Arsove},
Expo. Math. \textbf{7}(1989), 69--72.
\item[{[Ri99]}] J. Riihentaus, \emph{Subharmonic functions: non-tangential and 
tangential boundary
behavior},  Function Spaces, Differential Operators and Nonlinear Analysis,
Proceedings of the Sy\"ote Conference 1999 (eds.  V. Mustonen and J. R\'akosnik), 
 Math. Inst., Czech Acad. Science,  Praha, 2000, pp. 229--238.
\item[{[Ri01]}] J. Riihentaus, \emph{A generalized mean value inequality for 
subharmonic functions}, Expo. Math.  \textbf{19}(2001), 187--190.
\item[{[St98]}] M. Stoll, \emph{Weighted tangential boundary limits of subharmonic functions on
domains in ${\mathbb {R}}^n$ ($n\geq 2$)}, Math. Scand.  \textbf{83}(1998), 300--308.
\item[{[Su90]}]  N. Suzuki, \emph{Nonintegrability of harmonic functions in a domain},
Japan. J. Math. \textbf{16}(1990), 269--278.
\item[{[Su91]}] N. Suzuki, \emph{Nonintegrability of superharmonic functions}, 
Proc.\ Amer.\ Math.\ Soc. \ \textbf{113}(1991), 113--115.
\item[{[To86]}]  A. Torchinsky, \emph{Real-Variable Methods in Harmonic Analysis},
 Academic Press,  London, 1986.
\end{enumerate}
}

\vspace{1.0cm}
\footnotesize{
\noindent\textsc{South Carelia Polytechnic}\\
\textsc{P.O. Box 99}\\
\textsc{FIN-53101 Lappeenranta, Finland}\\

\noindent\textsc{and}\\

\noindent\textsc{Department of Mathematics}\\
\textsc{University of Joensuu}\\
\textsc{P.O. Box 111}\\
\textsc{FIN-80101 Joensuu, Finland}\\

\noindent E-mail address: \texttt{juhani.riihentaus{@}scp.fi}
}

\end{document}